\documentclass[11pt,reqno]{amsart}
\usepackage[left=30truemm,right=30truemm,top=30truemm,bottom=30truemm]{geometry}
\usepackage{eufrak}
\usepackage{amsfonts}
\usepackage{float}
\usepackage{amsmath,amssymb}
\usepackage{lscape}
\usepackage[dvipdfmx]{graphicx, color}
\usepackage[all]{xy}

\usepackage{comment}

\def\qed{\hfill $\Box$}
\def\proof{\noindent {\sl Proof} :\;  }

\newcommand{\Pcal}{\mathcal{P}}

\newcommand{\W}{\mathcal{W}}
\newcommand{\X}{\mathcal{X}}

\newcommand{\D}{\mathcal{D}}

\newcommand{\E}{\mathcal{E}}
\newcommand{\F}{\mathcal{F}}

\newcommand{\C}{\mathbb{C}}

\newcommand{\R}{\mathbb{R}}

\def\qed{\hfill $\Box$}
\def\proof{\noindent {\sl Proof} :\;  }

\def\rd{\partial}

\def\ba{\mbox{\boldmath $a$}}

\def\b0{\mbox{\boldmath $0$}}

\def\beeta{\mbox{\boldmath $\eta$}}
\def\bzeta{\mbox{\boldmath $\zeta$}}

\def\bba{\mbox{\tiny$\ba$}}

\def\bbeeta{\mbox{\tiny$\beeta$}}

\theoremstyle{plain}
\newtheorem{thm}{\bf Theorem}[section]

\newtheorem{lem}[thm]{\bf Lemma}
\newtheorem{prop}[thm]{\bf Proposition}

\theoremstyle{remark}
\newtheorem{dfn}[thm]{\bf Definition}

\newtheorem{rem}[thm]{\bf Remark}
\newtheorem{exam}[thm]{\bf Example}

\newcommand{\vp}{\varphi}

\newenvironment{talign}
 {\align}
 {\endalign}
\newenvironment{talign*}
 {\csname align*\endcsname}
 {\endalign}

 \makeatletter
 
 \@addtoreset{equation}{section}
\makeatother

\begin{document}
\title[The space of positive transition measures on a Markov chain]
{The space of positive transition measures on a Markov chain}
\author[N.~Nakajima]{Naomichi Nakajima}
\address[N.~Nakajima]{Department of Architecture, School of Architecture, Shibaura Institute of Technology, Tokyo, 135-8548, Japan}
\email{naomichi@shibaura-it.ac.jp}

\subjclass[2020]{Primary~53B12; Secondary~53A15}
\keywords{information geometry, Markov chains, transition probabilities, positive measures, Bregman divergences, $F$-divergences, the Perron-Frobenius theorem}
\dedicatory{}
\begin{abstract}
Information geometry of Markov chains has been studied using the dually flat structure of the space of transition probabilities. 
Although applications of this structure have been investigated, few attempts have examined its statistical meaning.
In this paper, we construct a foundation for investigating the statistical meaning based on Amari's theory of {\em positive measures}.
For the space of discrete distributions, Amari has introduced {\em the space of positive measures} by removing the constraint condition and investigated the extended space by finding the Bregman and $F$-divergence suitably. 
According to this, we introduce an extension of the space of transition probabilities equipped with suitable $F$-divergence for a given Markov chain.
We regard it as {\em the space of positive transition measures on a Markov chain}, and study its dually flat structure.
This provides new insight into the geometry of Markov chains and may lead to the development of the theory of Markov embeddings.
\end{abstract}

\maketitle


\section{Introduction}\label{sec:intro}
Information geometry of Markov chains has been investigated in several authors so far (\cite{Feigin, Hayashi2016, IA, KS1997, KS1998, Nagaoka, Takeuchi}).
Given a Markov chain, let $\W$ be the space of transition probabilities.
Then, a submanifold of $\W$ is called a Markov model.
The space $\W$ itself admits an exponential family, and thus it has the {\em dually flat structure} \cite{Nagaoka}.
In this framework, the dual potential function $\vp$ on the expectation parameter space $M$ of $\W$ takes a central role.

Although practical applications of this dually flat structure have been investigated (e.g., \cite{Hayashi2016}), its statistical meaning still seems to be unclear. 
For the space of probability distributions on a finite set, a well-known theorem due to Chentsov characterizes the Fisher-Rao metric and the Amari-Chentsov cubic tensor as the invariant structure of the space under Markov embeddings \cite{Chentsov}. 
This characterization essentially means that the dually flat structure is statistically natural. 
However, there have been few attempts to give such characterizations for Markov chains. 
While Markov embeddings on Markov models have been recently proposed by Wolfer and Watanabe based on purely statistical properties of Markov chains, this approach is still under development \cite{WW2024}.
In this paper, from a different viewpoint, we attempt to construct a foundation for investigating the statistical meaning of the dually flat structure of $\W$ based on Amari's well-established theory of {\em positive measures} on a finite set $S$ (cf. \cite{Amari2009, Amari16, AmariNagaoka00}).

The space $\Pcal(S)$ of probability distributions on $S$ is extended to the space of positive measures, denoted by $\bar{\Pcal}(S)$, by removing the constraint condition.
An $F$-divergence on $\bar{\Pcal}(S)$ is defined by using a given strictly convex function $F$ on $(0, \infty)$ with certain conditions, and it is invariant under Markov embeddings \cite{Amari2009}.
It is known that the dually flat structure of $\Pcal(S)$, which consists of the Fisher-Rao metric and the Amari-Chentsov cubic tensor, is naturally extended to that of $\bar{\Pcal}(S)$, and
the KL-divergence on $\bar{\Pcal}(S)$ and its restriction to $\Pcal(S)$ are Bregman divergences with respect to their dually flat structures \cite{Amari2009}.
Importantly, the KL-divergence is also an $F$-divergence, and in this sense, the KL-divergence is interpreted as a divergence which simultaneously derives the dually flat structures on both $\bar{\Pcal}(S)$ and $\Pcal(S)$ such that they are invariant under Markov embeddings.
This characterization using $F$-divergences is essentially the same as the one due to Chentsov.

We will develop the counterpart for the space $\W$ of transition probabilities on a Markov chain.
First, we extend $\W$ to a bigger space $\F^+$ which is obtained by removing the constraint condition for transition probabilities, and then  we define an $F$-divergence on $\F^+$.
Similarly, removing the conditions of the expectation parameter space $M$ derived from a given Markov chain being a probability distribution and stationary, we obtain the fully extended expectation parameter space $\overline{M}$, which is diffeomorphic to $\F^+$ (Lemma \ref{lem:T_bar}).
Then we give a divergence that is both a Bregman divergence and an $F$-divergence, which takes a similar role as the KL-divergence on $\bar{\Pcal}(S)$ (Theorem \ref{thm:Bregman}).
This divergence also restores the canonical divergence of $\W$ due to Nagaoka by restricting it to $\W$.
To show that the divergence is a Bregman divergence, we construct a (dual) potential function $\bar{\vp}$ on $\overline{M}$ explicitly.
Actually, the potential function has a $1$-dimensional kernel of its Hessian matrix at every point of $\overline{M}$, thus we take a hyperplane section $\tilde{M}$ in $\overline{M}$ so that a genuine dually flat structure is defined on it.
That induces a hypersurface $\tilde{\W}$ in $\F^+$ which should be a right object as the space of positive measures for the Markov chain.
To summarize this argument, we draw the following diagram:
$$\xymatrix{
\F^+ \ar@{}[d]|{\bigcup} \ar[r]_{\bar{T}}^{\sim} & \overline{M}  \ar@{}[d]|{\bigcup} \ar[rd]^{\bar{\vp}} & \\
\tilde{\W}  \ar[d] \ar[r]_{\bar{T}|_{\tilde{\W}}}^{\sim} & \tilde{M} \ar[d]  \ar[r]_{\bar{\vp}|_{\tilde{M}}} & \R \ar@{=}[d]\\
\W \ar[r]_{T}^{\sim} & M \ar[r]_{\vp} & \R 
}$$
The bottom row corresponds to the dually flat structure of Nagaoka \cite{Nagaoka} (i.e., the expectation parameter space $M$ and the dual potential function $\vp$), and our extension is the middle row.
Eventually, we claim that the pair $(\W, \tilde{\W})$ behaves like $(\Pcal(S), \bar{\Pcal}(S))$, and we call $\tilde{\W}$ the space of {\em positive transition measures}.

This paper is organized as follows.
In $\S 2$ we review the dually flat structure of a Markov model according to \cite{Nagaoka}.
In $\S 3$ we define the class of $F$-divergences on $\F^+$, and give a divergence which is both an $F$-divergence and a Bregman divergence.
Then the space of positive transition measures is introduced.
It is shown that the restriction of its dually flat structure to $\W$ restores that of $\W$ due to Nagaoka.
In $\S 4$, we discuss some statistical aspects of our theory.

\section{Transition probabilities on a Markov chain and the dually flat structure}\label{subsec:Markov}
In this section, we give our setup based on \cite{Nagaoka}. 
Let $\X := \{0, 1, \cdots, d\}$ $(d\geq 1)$ and $\E \subset \X\times\X$.
We consider Markov chains on the directed graph $(\X, \E)$, that is, we regard $\X$ as the state space, and the transition probabilities are defined on $\E$.
Let $\F^+$ denote the set of positive functions on $\E$ and $\W=\W(\X, \E)$ the set of transition probabilities:
$$
\F^+ = \{w:\E\to\R \mid w(x, y) > 0 \mbox{~for any~} (x, y) \in\E \},
$$
and
$$\textstyle
\W = \{w\in\F^+ \mid \sum_{y:(x,y)\in\E}w(x, y) = 1 \mbox{~for any~} x\in\X \}.
$$
Here, the notation 
$$
\sum_{y:(x, y) \in \E}
$$ 
means that fixing $x$, the sum runs over $y \in \X$ satisfying $(x, y) \in \E$  (so the sum depends on $x$). Also $\sum_{x:(x, y) \in \E}$ is similar. 
Throughout the present paper, we assume that $\E$ is strongly connected, that is, for any $x, y\in\X$ there exist $(x_1, x_2), (x_2, x_3), \cdots, (x_{N-1}, x_N) \in \E$ such that $x_1=x, x_N=y\;(N\geq 2)$.
This assumption means that for every $f\in\F^+$, the associated matrix $A(f)=[a_{ij}(f)]_{0\leq i,j\leq d}$ defined by 
$$
a_{ij}(f)=
\left\{
\begin{array}{ll}
f(i,j) & (i,j)\in\E\\
0 & (i,j)\notin\E
\end{array}
\right.
$$
is {\em irreducible} \cite{Zhan2013}. 
In particular, we call $A(w)$ a transition matrix for $w\in \W$. 
The following theorem is very important for studying properties of Markov chains.

\begin{thm}[the Perron-Frobenius theorem ,e.g., {\cite[Theorem 6.8]{Zhan2013}}]\label{PF}
Let $A$ be an $n\times n$ matrix. If all components of $A$ are non negative and $A$ is irreducible, then there exists a real eigenvalue $r>0$ of $A$ such that the following properties hold:
\begin{enumerate}\renewcommand{\labelenumi}{(\arabic{enumi})}
\item the geometric multiplicity and the algebraic multiplicity of $r$ are both one, and $r \geq |\lambda|$ for any eigenvalues $\lambda$,
\item there exists a unique left eigenvector $\mu=(\mu_1,\cdots,\mu_n)^T$ associated with $r$ such that $\mu_i>0$ for any $i$ and $\sum_{i=1}^n\mu_i=1$,
\item if there exists an eigenvector of $A$ whose all components are non negative, then it is an eigenvector of the eigenvalue $r$.
\end{enumerate}
We call $r$ the Perron-Frobenius root of $A$.
\end{thm}
For a positive function $f$ on $\E$, we write $r(f)$ and $\mu_f=(\mu_f(0),\cdots,\mu_f(d))^T$ as $r$ and $\mu$ above corresponding to $A(f)$, respectively. 
In this paper, abusing words, we call $\mu_f$ {\em the stationary distribution} for $f$ (not necessarily $f\in\W$).
We use the following lemma later.

\begin{lem}\label{lem:PF}
For $f\in\F^+$ and $a>0$, it holds that $r(af)=ar(f)$ and $\mu_{af}=\mu_f$.
\end{lem}
\proof
Since $\mu_f$ is the stationary distribution for $f$, it satisfies $\mu_f^TA(f)=r(f)\mu_f^T$.
By multiplying the both sides of this equation by $a$, we have $\mu_f^TA(af)=ar(f)\mu_f^T$.
Since all components of $\mu_f$ are positive, it follows that the eigenvalue for $\mu_f$ is the Perron-Frobenius root of $A(af)$ from Theorem \ref{PF} (3), i.e., $r(af)=ar(f)$.
Thus we get $\mu_f^TA(af)=r(af)\mu_f^T$, which means that $\mu_f$ is also the stationary distribution for $af$, i.e., $\mu_{af}=\mu_f$.
\qed

\begin{rem}
The Perron-Frobenius root $r(f)$ smoothly depends on $f\in\F^+$ and thus so does $\mu_f$. 
To verify this, let $P:\F^+\times(0,\infty)\to\R$ be the characteristic polynomial of $A(f)$, that is, $P(f, \lambda)=\det(\lambda I - A(f))$, where $I$ is the identity matrix.
Fix $f_0\in\F^+$ and put $r_0:=r(f_0)$. Since $r_0$ is a simple root of $P(f_0, \lambda)=0$, it holds that $P(f_0, r_0)=0$ and $\frac{\rd}{\rd\lambda}P(f_0, r_0)\neq 0$. 
From the implicit function theorem, there exist neighborhoods $U\subset\F^+$ of $f_0$, $V\subset(0,\infty)$ of $r_0$ and a smooth function $\tilde{r}:U\to V$ such that the set of zeros of $P(f,\lambda)$ in $U\times V$ coincides with the graph of $\tilde{r}$. 
On the other hand, from the continuity of roots of the polynomial $\det(\lambda I - A(f))$ (e.g., \cite{CC}), there exist continuous functions $z_0,\cdots,z_d:\F^+\to\C$ such that $\det(z_i(f) I - A(f))=0$ for any $f$. Since the Perron-Frobenius root $r(f)$ can be written as $r(f)=\max\{|z_0(f)|, \cdots, |z_d(f)|\}$, we see that $r:\F^+\to(0,\infty)$ is continuous. 
By replacing $U$ with smaller one if necessary, we can assume that $r(U)\subset V$. Obviously, $P(f, r(f))=0$ for $f\in U$. Hence, $r(f)=\tilde{r}(f)$ on $U$. 
\end{rem}

For each $w\in\W$, a joint probability distribution $p_w^{(n)}$ of $x^n:=(x_1,\cdots,x_n)$ satisfying $(x_i, x_{i+1})\in\E$ for all $i$ is defined by
$$
p_w^{(n)}(x^n) = \mu_w(x_1)w(x_1, x_2)\cdots w(x_{n-1},x_n),
$$
where $\mu_w=(\mu_w(0),\cdots,\mu_w(d))^T$ is the stationary distribution for $w$.
Nagaoka \cite{Nagaoka} focused on the space of conditional probabilities $\W$, not on that of probability distributions $\{p_w^{(n)}\}$, which is depending on the sample size $n$.
Now a family $\{w_\theta\in\W\}_{\theta\in\R^k}$ of probability distributions is defined by 
$$
\log w_\theta(x, y) = C(x, y) + \sum_{i=1}^k\theta_iF_i(x, y) + K_\theta(y) - K_\theta(x) - \psi(\theta),
$$
where $\theta=(\theta_1,\cdots,\theta_k)\in\R^k$ are parameters (called {\em natural parameters}), $C, F_1,\cdots,F_k$ are functions on $\E$, $K:\X\times\R^k\to\R,\;(x,\theta)\mapsto K_\theta(x)$ and $\psi:\R^k\to\R$ is the normalization factor; it is called a {\em $k$-dimensional exponential family} of Markov chains on $(\X, \E)$ \cite{Nagaoka}.
This family induces the {\em dually flat structure} $(g, \nabla^{(e)}, \nabla^{(m)})$ of $\W$ as follows.
See \cite{Nagaoka} for the detail.
The metric $g=[g_{ij}(\theta)]_{1\leq i,j\leq k}$ is defined by
$$
g_{ij}(\theta):=\sum_{(x,y)\in\E}p_{w_\theta}^{(2)}(x,y)\left(\frac{\rd}{\rd\theta_i}\log w_\theta(x,y)\right)\left(\frac{\rd}{\rd\theta_j}\log w_\theta(x,y)\right),
$$
called the {\em Fisher metric}.
Indeed, it is the limit of the Fisher matrix $[g_{ij}^{(n)}(\theta)]$ of the family $\{p^{(n)}_{w_\theta}\}_{\theta\in\R^k}$ in the usual sense:
\begin{align*}
g_{ij}(\theta) = \lim_{n \to \infty}\frac{1}{n}g_{ij}^{(n)}(\theta).
\end{align*}
Let $\F$ denote the set of functions on $\E$ and put
\begin{align*}
\F^S&=\F^S(\X, \E):=\{\textstyle f\in\F \mid \sum_{y:(x,y)\in\E}f(x, y) = \sum_{y:(y,x)\in\E}f(y, x) \mbox{~for any~} x\in\X\},\\
\F^A&=\F^A(\X, \E):=\{f\in\F \mid f(x,y)=\kappa(y) - \kappa(x),\;\; \kappa:\X\to\R\}.
\end{align*}
Then, by letting $\R$ denote the space of constant functions on $\E$, it is shown that the quotient vector space $\F / (\F^A\oplus\R)$ is diffeomorphic to $\W$ and the affine structure of $\F / (\F^A\oplus\R)$ gives the natural parameters $\theta$ of $\W$ as an exponential family.
Also, the mapping 
\begin{align*} 
T: \W\to\Pcal(\E)\cap\F^S, \;\;w\mapsto p^{(2)}_w = (\mu_w(x)w(x,y))_{(x,y)\in\E}
\end{align*}
with $\Pcal(\E)=\{f\in\F^+\mid \sum_{(x,y)\in\E}f(x,y)=1\}$ is a diffeomorphism, and it gives the new coordinates of $\W$, called {\em the expectation parameters}.
The natural parameters and the expectation parameters are affine coordinates with respect to the affine connections $\nabla^{(e)}$ and $\nabla^{(m)}$ defined by
\begin{align*}
\Gamma^{(e)}_{ij, k} 
&:= g(\nabla^{(e)}_{\rd_i}\rd_j, \rd_k) 
= \sum_{(x, y)\in\E}\rd_i\rd_j\log w_\theta(x, y)\rd_kp^{(2)}_{w_\theta}(x, y),\\
\Gamma^{(m)}_{ij, k} 
&:= g(\nabla^{(m)}_{\rd_i}\rd_j, \rd_k)
= \sum_{(x, y)\in\E}\rd_i\rd_jp^{(2)}_{w_\theta}(x, y)\rd_k\log w_\theta(x, y),
\end{align*}
where we put $\rd_i:=\frac{\rd}{\rd \theta_i}$ for simplicity. These connections are dual to each other with respect to the Fisher metric $g$, that is, it holds that 
$$
Xg(Y, Z)=g(\nabla^{(e)}_XY, Z) + g(Y, \nabla^{(m)}_XZ)
$$
for any vector fields $X, Y, Z$ on $\W$. 
For the dually flat structure $(g, \nabla^{(e)}, \nabla^{(m)})$, the (dual) Bregman divergence $\D:\W\times\W\to\R$ is defined by
\begin{align}\label{divergence}
\D(w_1, w_2) 
= \sum_{(x, y)\in\E} \mu_{w_1}(x)w_1(x,y)\log\frac{w_1(x, y)}{w_2(x, y)}.
\end{align}

\section{Positive transition measures on a Markov chain and the dually flat structure}\label{sec:positive}
For a finite set $S=\{0,\cdots,n\}$, let $\Pcal(S)$ and $\bar{\Pcal}(S)$ denote the space of probability distributions and the space of positive measures on $S$, respectively:
\begin{align*}
\Pcal(S)&=\{(p_0,\cdots,p_n)\in\R^{n+1}\mid p_i>0,\;\;\sum_{i=0}^n p_i=1\},\\
\bar{\Pcal}(S)&=\{(p_0,\cdots,p_n)\in\R^{n+1}\mid p_i>0\}.
\end{align*}
Given a strictly convex function $F:(0, \infty)\to\R$ with $F(1)=F'(1)=0$ and $F''(1)=1$, called {\em a standard convex function} (\cite{Amari2009}), the function $\D_F:\bar{\Pcal}(S)\times\bar{\Pcal}(S)\to\R$ defined by
$$
\D_F(p, q)=\sum_{i=0}^np_iF\left(\frac{q_i}{p_i}\right)
$$
is called {the $F$-divergence} on $\bar{\Pcal}(S)$, where $p=(p_0,\cdots,p_n), q=(q_0,\cdots,q_n)$. In the case where 
$$
F(t)=-\log t + (t-1),
$$
the $F$-divergence $\D_F$ is the KL-divergence: 
$$
\D_F(p,q)=\sum_{i=0}^np_i\log\left(\frac{p_i}{q_i}\right) + \sum_{i=0}^nq_i - \sum_{i=0}^np_i.
$$
Amari has characterized the dually flat structures of $\bar{\Pcal}(S)$ and $\Pcal(S)$ by finding the Bregman and $F$-divergences.
In fact, the KL-divergence and its restriction to $\Pcal(S)$ are also Bregman divergences with respect to their dually flat structures which are invariant under Markov embeddings \cite{Amari2009}.
The readers are referred to \cite{Amari2009} for more detail.

We attempt to develop a foundation for studying such characterizations of $\W$ by using the bigger space $\F^+$.

\begin{dfn}\label{def:f-div}
Let $F:(0,\infty)\to\R$ be a standard convex function. We define the $F$-divergence on $\F^+$ as $\D_F:\F^+\times\F^+\to\R$,
\begin{align*}
\D_F(f, g) = \sum_{(x,y)\in\E}\mu_f(x)f(x, y)F\left( \left.\frac{g(x, y)}{r(g)} \right/ \frac{f(x, y)}{r(f)}\right).
\end{align*}
\end{dfn}
\noindent
In the context of {\em statistical manifolds}, it is well-known that a symmetric $(0, 2)$-tensor on a smooth manifold $N$ is induced by a some `asymmetric distance function' $\rho:N\times N\to\R$ as follows. 
See \cite{Eguchi} for the detail.
For a function $\rho$ and vector fields $X_1, \cdots, X_k, Y_1,\cdots, Y_l$ on $N$, we set a function 
$$
\rho[X_1\cdots X_k|Y_1\cdots Y_l]:N\to\R
$$
defined by
\begin{align*} 
\rho[X_1\cdots X_k|Y_1\cdots Y_l](r) = (X_1)_p\cdots(X_k)_p(Y_1)_q\cdots(Y_l)_q(\rho(p,q))|_{p=q=r}.
\end{align*}
In particular, $\rho$ is called a {\em weak contrast function} on $N$ (\cite{NO}) if it holds that for each $r\in N$
\begin{enumerate}
\item[--] $\rho[-|-](r)=\rho(r,r) = 0$ and 
\item[--] $\rho[X|-](r) = \rho[-|X](r) = 0$.
\end{enumerate}
Then a symmetric $(0, 2)$-tensor $h$ on $N$ is defined by
$$
h(X, Y):=-\rho[X|Y]=\rho[XY|-]=\rho[-|XY].
$$

\begin{prop}\label{prop:F-div}
The $F$-divergence $\D_F$ has the following properties:
\begin{enumerate}\renewcommand{\labelenumi}{(\arabic{enumi})}
\item $\D_F(f,g)\geq 0$.
\item $\D_F(f, g)=0$ if and only if $g=af$ for some $a>0$.
\item $\D_F$ is a weak contrast function on $\F^+$. Let $h_F$ denote the symmetric $(0, 2)$-tensor on $\F^+$ induced by $\D_F$.
\item The null space of $h_F$ at $f\in\F^+$ is the tangent space of the half line $\{af\mid a>0\}\subset\F^+$.
\end{enumerate}
\end{prop}
\proof
The non-negativity of $F$ yields (1).
We show (2).
Note that $r(af)=ar(f)$ for $a>0$ from Lemma \ref{lem:PF}.
Obviously if $g=af$ for some $a>0$, then $\D_F(f,g)=0$. 
Conversely, assume that $\D_F(f, g)=0$. Then $F(\frac{g(x,y)}{r(g)}/\frac{f(x,y)}{r(f)})=0$ for each $(x,y)\in\E$.
Since $F$ is a standard convex function, $F(t)=0$ if and only if $t=1$. Hence we get $g(x,y)=\frac{r(g)}{r(f)}\cdot f(x,y)$ for each $(x,y)\in\E$.
Therefore we see that (2) holds.
From (1) and (2), $\D_F$ takes the minimum value zero on the diagonal set in $\F^+\times\F^+$, and thus we get for each $f\in\F^+$
\begin{enumerate}
\item[--] $\D_F(f, f)=0$ and 
\item[--] $\D_F[X|-](f) = \D_F[-|X](f) = 0$,
\end{enumerate}
where $X$ is a vector field on $\F^+$.
This means that (3) is true.
Finally we show (4).
We consider a natural system of coordinates $\ba=(a_{xy})_{(x,y)\in\E}$ of $\F^+$ defined by $a_{xy}(f)=f(x,y)$.
Fix a point $f\in\F^+$.
By differentiating $\D_F(f, \cdot):\F^+\to\R$ at $f$, we have
\begin{align*}
&\D_F[-|\tfrac{\rd}{\rd a_{st}}](f)
= \left.\left\{\sum_{(x,y)\in\E}\mu_f(x)f(x,y)F'\left( \left.\tfrac{a_{xy}}{r(\bba)} \right/ \tfrac{f(x,y)}{r(f)}\right)\cdot\tfrac{\rd}{\rd a_{st}}\left(\tfrac{a_{xy}}{r(\bba)}\right)\right\}\right|_{\bba=f}=0,\\
&\D_F[-|\tfrac{\rd}{\rd a_{uv}}\tfrac{\rd}{\rd a_{st}}](f)
= \sum_{(x,y)\in\E}\mu_f(x)f(x, y)\left.\tfrac{\rd}{\rd a_{st}}\left(\tfrac{a_{xy}}{r(\bba)}\right)\right|_{\bba=f}\cdot\left.\tfrac{\rd}{\rd a_{uv}}\left(\tfrac{a_{xy}}{r(\bba)}\right)\right|_{\bba=f}
\end{align*}
for $F'(1)=0$ and $F''(1)=1$.
Hence for $X=\sum_{(x,y)\in\E}v_{xy}(\frac{\rd}{\rd a_{xy}})_{f}\in T_f\F^+$ we get 
$$
h_F(X, X) = \sum_{(x,y)\in\E}\mu_f(x)f(x,y)\left\{\sum_{(s,t)\in\E}\left.v_{st}\tfrac{\rd}{\rd a_{st}}\left(\tfrac{a_{xy}}{r(\bba)}\right)\right|_{\bba=f}\right\}^2.
$$
Since
$$
\left.\frac{\rd}{\rd a_{st}}\left(\frac{a_{xy}}{r(\ba)}\right)\right|_{\bba=f}
= \frac{\delta_{sx}\delta_{ty}r(f)-f(x,y)\cdot\frac{\rd r}{\rd a_{st}}(f)}{(r(f))^2},
$$
we see that $h_F(X, X)=0$ if and only if it holds that for each $(x,y)\in\E$
\begin{align*}
\sum_{(s,t)\in\E}\left.v_{st}\tfrac{\rd}{\rd a_{st}}\left(\tfrac{a_{xy}}{r(\bba)}\right)\right|_{\bba=f} = \frac{v_{xy}}{r(f)} - \frac{f(x,y)}{(r(f))^2}\sum_{(s,t)\in\E}v_{st}\frac{\rd r}{\rd a_{st}}(f) = 0.
\end{align*}
Note that $\sum_{(s,t)\in\E}v_{st}\frac{\rd r}{\rd a_{st}}(f)$ is independent of $(x,y)\in\E$, thus it is constant.
Hence the vector $X=\sum_{(x,y)\in\E}v_{xy}(\frac{\rd}{\rd a_{xy}})_{f}$ is parallel to $\sum_{(x,y)\in\E}f(x,y)(\frac{\rd}{\rd a_{xy}})_{f}$.
We have thus proved the proposition.
\qed

In order to give an $F$-divergence on $\F^+$ which is also a Bregman divergence, we consider the extended space of expectation parameters
$$
\overline{M}:=\{\beeta=(\eta_{xy})_{(x,y)\in\E} \mid \eta_{xy}>0\}.
$$
The genuine expectation parameter space $M$ is the affine subspace of $\overline{M}$ defined by the two conditions:
\begin{talign}\textstyle
&\sum_{(x,y)\in\E}\eta_{xy}=1 \mbox{~and~}\label{normalization condition}\\
&\sum_{y:(x,y)\in\E}\eta_{xy}=\sum_{y:(y,x)\in\E}\eta_{yx} \mbox{~for any~} x\in\X.\nonumber
\end{talign}
Note that $M=\Pcal(\E)\cap\F^S$.
In particular, the first condition is derived from $\Pcal(\E)$, and we call it {\em the normalization condition} on $\overline{M}$.
We set 
$$
\bar{T}:\F^+\to\overline{M},\;\; f\mapsto (\mu_f(x)f(x,y))_{(x,y)\in\E}.
$$
We also set for $\beeta=(\eta_{xy})_{(x,y)\in\E}\in\overline{M}$
$$
r(\beeta):=\sum_{(x,y)\in\E}\eta_{xy},\quad
\eta^x:=\sum_{k:(k, x)\in\E}\eta_{kx} \;\;\mbox{~and~}\;\;
\eta_x:=\sum_{k:(x, k)\in\E}\eta_{xk}.
$$

\begin{lem}\label{lem:T_bar}
$\bar{T}$ has the following properties:
\begin{enumerate}\renewcommand{\labelenumi}{(\arabic{enumi})}
\item $\bar{T}|_{\W}=T:\W\overset{\sim}{\to} M$,
\item $\bar{T}$ is a diffeomorphism; its inverse, denoted by $\bar{\tau}:\overline{M}\to\F^+$, is given by 
\begin{align*}
\bar{\tau}(\beeta):\E\to\R,\;\;(x, y)\mapsto r(\beeta)\frac{\eta_{xy}}{\eta^x}. 
\end{align*}
\item $\bar{T}(af)=a\bar{T}(f)$ for $f\in\F^+$ and $a>0$.
\end{enumerate}
\end{lem}
\proof
From the form of the mapping $\bar{T}$, the equality (1) is obvious.
Also, Lemma \ref{lem:PF} shows the formula (3).
We show (2) below. 
Take an arbitrary $f\in\F^+$ and put $\beeta:=\bar{T}(f)=(\eta_{xy})$ with $\eta_{xy}=\mu_f(x)f(x,y)$.
Then for each $x\in\X$
$$
\eta^x = \sum_{k:(k,x)\in\E}\eta_{kx} = \sum_{k:(k,x)\in\E}\mu_f(k)f(k, x) = r(f)\mu_f(x),
$$
where $r(f)$ is the Perron-Frobenius root for $f$.
Then, the vector $(\eta^x)_{x\in\X}$ is a left eigenvector of $A(f)$ for $r(f)$, and 
$$
\sum_{x\in\X}\eta^x = \sum_{(x,y)\in\E}\eta_{xy} = r(\beeta).
$$
From Theorem \ref{PF} (3), the vector $(\frac{\eta^x}{r(\bbeeta)})_{x\in\X}$ coincides with $\mu_f=(\frac{\eta^x}{r(\bbeeta)})_{x\in\X}$, that is, $r(\beeta)=r(f)$.
Hence, $f(x,y) = \frac{\eta_{xy}}{\mu_f(x)}=r(\beeta)\frac{\eta_{xy}}{\eta^x}$, thus $\bar{\tau}(\beeta)=f$. It is easily seen that $\bar{T}\circ\bar{\tau}(\beeta)=\beeta$.
Thus (2) is proven.
\qed

Here, we summarize the relations between $f\in\F^+$ and $\beeta=\bar{T}(f)\in\overline{M}$, which are obtained in the proof above:
\begin{align}\label{relations}
r(\beeta)=r(f),\quad
\eta^x= r(f)\mu_f(x),\quad
\eta_{xy}=\mu_f(x)f(x, y).
\end{align}

\begin{thm}\label{thm:Bregman}
Let $F(t)=-\log t + (t-1)$. Then the $F$-divergence is the Bregman divergence given by the following potential function on $\overline{M}$:
\begin{align*}
\bar{\vp}(\beeta) = \sum_{(x,y)\in\E}\eta_{xy}\log\eta_{xy} - \sum_{x\in\X}\eta_x\log\eta^x.
\end{align*}
\end{thm}

\proof
This is shown by direct computations. 
The $F$-divergence $\D_F$ is written as follows:
\begin{align*}
\D_F(f, g) 
= &\sum_{(x,y)\in\E}\mu_f(x)f(x, y)F\left( \left.\frac{g(x, y)}{r(g)} \right/ \frac{f(x, y)}{r(f)}\right)\\
= &\sum_{(x,y)\in\E}\mu_f(x)f(x, y)\left\{\log\left(\frac{r(g)f(x,y)}{r(f)g(x,y)}\right) + \frac{r(f)g(x,y)}{r(g)f(x,y)}-1   \right\}\\
= &\sum_{(x,y)\in\E}\mu_f(x)f(x, y)\log\frac{f(x,y)}{g(x,y)} - r(f)\log\frac{r(f)}{r(g)} \\
&\quad+ \frac{r(f)}{r(g)}\sum_{(x,y)\in\E}\mu_f(x)g(x, y) - r(f).
\end{align*}
On the other hand, for $\beeta=(\eta_{xy})=\bar{T}(f)$ and $\bzeta=(\zeta_{xy})=\bar{T}(g)$, the Bregman divergence $\D_{\rm{Bre}}:\overline{M}\times\overline{M}\to\R$ is given by
\begin{align*}
\D_{\rm{Bre}}(\beeta, \bzeta) 
&= \bar{\vp}(\beeta) - \bar{\vp}(\bzeta) + \sum_{(x,y)\in\E}\frac{\rd\bar{\vp}}{\rd\eta_{xy}}(\bzeta)(\zeta_{xy}-\eta_{xy})\\
&= \sum_{(x,y)\in\E}\eta_{xy}\log\eta_{xy} - \sum_{x\in\X}\eta_x\log\eta^x
-\sum_{(x,y)\in\E}\zeta_{xy}\log\zeta_{xy} + \sum_{x\in\X}\zeta_x\log\zeta^x\\
&\qquad+\sum_{(x,y)\in\E}\left(\log\zeta_{xy} -\log\zeta^x - \frac{\zeta_y}{\zeta^y} + 1\right)(\zeta_{xy}-\eta_{xy})\\
&= \sum_{(x,y)\in\E}\eta_{xy}\log\frac{\eta_{xy}}{\zeta_{xy}} - \sum_{x\in\X}\eta_x\log\frac{\eta^x}{\zeta^x} + \sum_{x\in\X}\eta^x\frac{\zeta_x}{\zeta^x} - r(\beeta).
\end{align*}
For the last equality above, we used some formulas, for example,
$$
\sum_{(x,y)\in\E}\eta_{xy}\log\zeta^x = \sum_{x\in\X}\left(\sum_{y:(x,y)\in\E}\eta_{xy}\right)\log\zeta^x = \sum_{x\in\X}\eta_x\log\zeta^x
$$
and
$$
\sum_{(x,y)\in\E}\eta_{xy}\frac{\zeta_y}{\zeta^y}=\sum_{y\in\X}\left(\sum_{x:(x,y)\in\E}\eta_{xy}\right)\frac{\zeta_y}{\zeta^y} = \sum_{x\in\X}\eta^x\frac{\zeta_x}{\zeta^x}.
$$
Using the relations (\ref{relations}), we have 
\begin{align*}
\D_{\rm{Bre}}(\bar{T}(f), \bar{T}(g))
&= \sum_{(x,y)\in\E}\eta_{xy}\log\frac{\mu_f(x)f(x,y)}{\mu_g(x)g(x,y)}
- \sum_{x\in\X}\eta_x\log\frac{\mu_f(x)r(f)}{\mu_g(x)r(g)}\\
&\quad\,+\sum_{x\in\X}\mu_f(x)r(f)\frac{\zeta_x}{\mu_g(x)r(g)} - r(f)\\
&\quad= \sum_{(x,y)\in\E}\eta_{xy}\log\frac{f(x,y)}{g(x,y)}
- r(f)\log\frac{r(f)}{r(g)}\\
&\qquad\,+ \frac{r(f)}{r(g)}\sum_{x\in\X}\left(\sum_{y:(x,y)\in\E}\mu_g(x)g(x, y)\right)\frac{\mu_f(x)}{\mu_g(x)} - r(f)\\
&= \sum_{(x,y)\in\E}\mu_f(x)f(x,y)\log\frac{f(x,y)}{g(x,y)} - r(f)\log\frac{r(f)}{r(g)}\\
&\qquad\,+ \frac{r(f)}{r(g)}\sum_{(x,y)\in\E}\mu_f(x)g(x, y) - r(f)\\
&= \D_F(f, g).
\end{align*}
This completes the proof.\qed

\begin{rem}
The restriction of the $F$-divergence with $F(t)=-\log t + (t-1)$ to $M$ restores the divergence (\ref{divergence}) associated to the dually flat structure of Nagaoka. 
In fact, for $f, g\in\W$ we see 
\begin{align*}
\D_F(f,g)&=\sum_{(x,y)\in\E}\mu_f(x)f(x,y)F\left(\frac{g(x,y)}{f(x,y)}\right)\\
&= \sum_{(x,y)\in\E}\mu_f(x)f(x,y)\log\frac{f(x,y)}{g(x,y)}+\sum_{x\in\X}\sum_{y:(x,y)\in\E}\mu_f(x)(g(x,y)-f(x,y))\\
&= \sum_{(x,y)\in\E}\mu_f(x)f(x,y)\log\frac{f(x,y)}{g(x,y)}.
\end{align*}
\end{rem}

\begin{rem}\label{rem:Takeuchi}
We note that \cite{Takeuchi}, with a different motivation from ours, implicitly treats $\overline{M}$ as a parameter space of $\W$ (not $\F^+$) and introduces a symmetric $(0, 2)$-tensor $\bar{g}$ on $\overline{M}$ by pulling back $g$ on $M$ via a canonical linear projection $\overline{M}\to M$.
Thus $\bar{g}$ is degenerate along the kernel directions of the projection.
Also, \cite{Konno2023} gives the following potential function for $\bar{g}$, which is similar to ours:
\begin{align*}
\hat{\vp}(\beeta) = \sum_{(x,y)\in\E}\eta_{xy}\log\eta_{xy} - \sum_{x\in\X}\eta_x\log\eta_x,\;\;\beeta=(\eta_{xy})_{(x,y)\in\E}\in\overline{M}.
\end{align*}
\end{rem}

From a simple computation, we see that $\bar{\vp}$ is homogeneous of degree $1$, i.e., $\bar{\vp}(a\beeta)=a\bar{\vp}(\beeta)$ for all $a>0$ and $\beeta\in\overline{M}$.
Then we see that the Hessian matrix of $\bar{\vp}$ at every point $\beeta\in\overline{M}$ has the $1$-dimensional kernel spanned by the numerical vector $\beeta\in\R^{|\E|}\cong T_{\bbeeta}\overline{M}$.
Indeed, it follows from Lemma \ref{lem:T_bar} (3) that the $1$-dimensional null space of $h_F$ with $F(t)=-\log t + (t-1)$ at $f=\bar{T}^{-1}(\beeta)$ is sent to the kernel above by $d\bar{T}_f:T_f\F^+\overset{\sim}{\to}T_{\bbeeta}\overline{M}$.
Therefore, by imposing only the normalization condition (\ref{normalization condition}) on $\overline{M}$, we have the hyperplane section $\tilde{M}$ in $\overline{M}$ so that $\bar{\vp}$ is strictly convex on it:
$$
\tilde{M}:=\{\beeta=(\eta_{xy})\in\overline{M}\mid r(\beeta)=1\}.
$$
Using the relation $r(f)=r(\beeta)$ with $\bar{T}(f)=\beeta$ , we get the genuine dually flat manifold $\tilde{\W}$, which is an extended space of $\W$ as a hypersurface in $\F^+$:
$$
\tilde{\W}:=\{f\in\F^+\mid r(f)=1\}.
$$

\begin{thm}\label{cor:duallyflat}
The hypersurface $\tilde{\W}$ has the dually flat structure induced by the potential function $\tilde{\vp}:=\bar{\vp}|_{\tilde{M}}$ on $\tilde{M}$; the restriction of this dually flat structure to $\W$ restores the dually flat structure of Nagaoka.
We call $\tilde{\W}$ the space of positive transition measures on $(\X, \E)$.
Moreover $F$-divergences on $\tilde{\W}$ are written as
$$
\D_F(f, g) = \sum_{(x,y)\in\E}\mu_f(x)f(x,y)F\left(\frac{g(x,y)}{f(x,y)}\right),
$$
where $f,g\in\tilde{\W}$.
\end{thm}

\begin{exam}
Let us consider the case where $\X=\{0, 1\}$ and $\E=\X\times\X$.
We identify $\F^+$ with the open domain of $\R^4$ where all coordinates $(x,y,z,w)$ are positive by putting $x=f(0,0),\; y=f(0,1),\; z=f(1,0)$ and $w=f(1,1)$ for $f\in\F^+$.
Then the Perron-Frobenius root $r(f)$ for $f=(x,y,z,w)$ is given by
$$
r(f) = \frac{x+w+\sqrt{(x-w)^2+4yz}}{2}.
$$
The equation $r(f)=1$ is equivalent to
$$
x + w -2 = -\sqrt{(x-w)^2+4yz},
$$
which yields $x+w<2$.
Moreover, squaring both sides of the equation above we get
$$
(x-1)(w-1) - yz = 0.
$$
Therefore the space $\tilde{\W}$ of positive transition measures on $(\X, \E)$ is given by
$$
\tilde{\W}=\{f=(x,y,z,w)\in\R^4\mid (x-1)(w-1) - yz = 0,\; x+w<2 \mbox{~and~} x,y,z,w > 0\}.
$$
Also, we compute the Hessian matrix of $\tilde{\vp}$.
The second order partial derivatives of $\bar{\vp}$ at $\beeta=(\eta_{00},\; \eta_{01},\; \eta_{10},\; \eta_{11})\in \overline{M}$ are given by
\begin{align*}
\frac{\rd\bar{\vp}}{\rd\eta_{st}\rd\eta_{uv}}(\beeta) = \frac{\delta_{su}\delta_{tv}}{\eta_{st}} - \frac{\delta_{sv}}{\eta^s} - \frac{\delta_{tu}\eta^t-\delta_{tv}\eta_t}{(\eta^t)^2}.
\end{align*}
Thus the Hessian matrix of $\bar{\vp}$ is obtained by
\begin{align*}
\left[
\begin{array}{cccc}
\frac{1}{\eta_{00}} - \frac{1}{\eta^{0}} - \frac{\eta^{0} - \eta_{0}}{\left(\eta^{0}\right)^{2}} & - \frac{1}{\eta^{0}} & - \frac{1}{\eta^{0}} + \frac{\eta_{0}}{\left(\eta^{0}\right)^{2}} & 0\\
- \frac{1}{\eta^{0}} & \frac{1}{\eta_{01}} + \frac{\eta_{1}}{\left(\eta^{1}\right)^{2}} & - \frac{1}{\eta^{1}} - \frac{1}{\eta^{0}} & - \frac{1}{\eta^{1}} + \frac{\eta_{1}}{\left(\eta^{1}\right)^{2}}\\
- \frac{1}{\eta^{0}} + \frac{\eta_{0}}{\left(\eta^{0}\right)^{2}} & - \frac{1}{\eta^{1}} - \frac{1}{\eta^{0}} & \frac{1}{\eta_{10}} + \frac{\eta_{0}}{\left(\eta^{0}\right)^{2}} & - \frac{1}{\eta^{1}}\\
0 & - \frac{1}{\eta^{1}} + \frac{\eta_{1}}{\left(\eta^{1}\right)^{2}} & - \frac{1}{\eta^{1}} & \frac{1}{\eta_{11}} - \frac{1}{\eta^{1}} - \frac{\eta^{1} - \eta_{1}}{\left(\eta^{1}\right)^{2}}
\end{array}
\right],
\end{align*}
which has the kernel spanned by the vector $(\eta_{00},\; \eta_{01},\; \eta_{10},\; \eta_{11})^T$.
Taking a coordinate system $(\eta_{00},\; \eta_{01},\; \eta_{10})$ of $\tilde{M}\;(\eta_{11}=1-\eta_{00}-\eta_{01}-\eta_{10})$ and restricting the matrix above to $\tilde{M}$ we get the Hessian matrix of $\tilde{\vp}$ at $\beeta\in\tilde{M}$: 
\begin{align*}
\left[
\begin{array}{ccc}
\frac{1}{\eta_{11}} + \frac{1}{\eta_{00}} - \frac{2}{\eta^{0}\eta^{1}} + \frac{\eta_{1}}{\left(\eta^{1}\right)^{2}} + \frac{\eta_{0}}{\left(\eta^{0}\right)^{2}} & \frac{1}{\eta_{11}} - \frac{1}{\eta^{0}\eta^{1}} & \frac{1}{\eta_{11}} - \frac{1}{\eta^{0}\eta^{1}} + \frac{\eta_{1}}{\left(\eta^{1}\right)^{2}}  + \frac{\eta_{0}}{\left(\eta^{0}\right)^{2}}\\
\frac{1}{\eta_{11}} - \frac{1}{\eta^{0}\eta^{1}} & \frac{1}{\eta_{01}} + \frac{1}{\eta_{11}} & \frac{1}{\eta_{11}} - \frac{1}{\eta^{0}\eta^{1}}\\
\frac{1}{\eta_{11}} - \frac{1}{\eta^{0}\eta^{1}} + \frac{\eta_{1}}{\left(\eta^{1}\right)^{2}} + \frac{\eta_{0}}{\left(\eta^{0}\right)^{2}} & \frac{1}{\eta_{11}} - \frac{1}{\eta^{0}\eta^{1}} & \frac{1}{\eta_{11}} + \frac{1}{\eta_{10}} + \frac{\eta_{1}}{\left(\eta^{1}\right)^{2}} + \frac{\eta_{0}}{\left(\eta^{0}\right)^{2}}
\end{array}
\right].
\end{align*}
\end{exam}

\section{Discussions}
In this paper, given a Markov chain $(\X, \E)$, we considered the space $\F^+$ of all positive real-valued functions on $\E$, as the largest extension of the space $\W$ of transition probabilities.
We first defined the class of $F$-divergences on $\F^+$.
Then we gave such an $F$-divergence which is also a Bregman divergence by regarding $\overline{M}$ as the expectation parameter space of $\F^+$.
Moreover, we gave a dually flat manifold $\tilde{\W}$ which is an extension of $\W$ by analyzing the kernels of the potential function $\bar{\vp}$ on $\overline{M}$.
In what follows, we briefly discuss possible applications of our geometric framework to statistics of Markov chains.

\subsection{Characterization of the dually flat structure}
In order to establish our theory, we need to further discuss the statistics of Markov chains, as the statistics of discrete distributions is essential to Amari's theory.
In other words, the theory of positive transition measures will develop in relation to the statistics of Markov chains, such as sufficient statistics and Markov embeddings (cf. \cite{WW2024}).
In fact, our $F$-divergences should be characterized by information monotonicity and by Markov embeddings on Markov models (cf. \cite{AmariNagaoka00, Csiszar1977, Csiszar2008, Jiao2014}). 
Conversely, these divergences may provide new insights into the theory of Markov embeddings.
In this way, the reciprocal development between the theory of positive transition measures and that of Markov embeddings leads to the establishment of the statistical meaning of the dually flat structure of Markov models.
Thus, the information geometry of Markov chains that is consistent with statistics should be established.

\subsection{Estimation and testing based on $F$-divergences}

Even though we focused on the hypersurface $\tilde{\W}$ to seek an analogy to Amari's theory in this paper, the largest space $\F^+$ may also be useful for statistical estimation of transition probabilities.

Although $\F^+$ has no dually flat structure adapted to the potential function $\bar{\vp}$, it has a {\em quasi-Hessian structure}, which was recently introduced as a singular version of a dually flat structure in \cite{NO}.
Even on this singular structure, the generalized Pythagorean theorem and the projection theorem still hold, and these theorems lead to new methods of estimation and testing for transition probabilities.

In fact, it is known that the {\em toric homogeneous Markov chain model} (the THMC model, for short) is useful for goodness-of-fit tests, and it is parameterized by $\F^+$; see \cite{Takemura2012, KS1996} for the definition and detailed properties.
Based on our geometric structure, new estimation and testing methods for the THMC model are suggested as follows.

Rao introduced test statistics using the Fisher--Rao Riemannian distance to investigate coordinate-free properties of testing \cite{Rao}.
Analogously, our framework introduces test statistics using $F$-divergences and reformulates goodness-of-fit tests for the THMC model.  
Furthermore, our Bregman divergence enables us to estimate transition probabilities of Markov chains from positive transition measures of the THMC model.
In the extended expectation parameter space $\overline{M}$, this estimation is interpreted as the orthogonal projection to $M$ along straight lines.
These testing and estimation methods should have good statistical properties.

Indeed, \cite{Hayashi2016} introduced an estimation method based on the dually flat structure on $\W$ of \cite{Nagaoka}, and the resulting estimator is asymptotically efficient and has a lower computational cost than existing MLE procedures.
Moreover, our $F$-divergences are projective in the sense of Proposition \ref{prop:F-div}(2).
Projective divergences are effective for robust estimation---for example, the $\gamma$-divergence \cite{Fujisawa}, the pseudo-spherical divergence \cite{GR2007} and its dual \cite{Hino2023}; hence our $F$-divergences may be robust for estimating transition probabilities.
From the viewpoint of information geometry, robust estimation should be investigated via the quasi-Hessian structures induced by projective divergences, and this investigation will clarify the statistical properties of our methods.

\subsection*{Acknowledgement}
The author would like to thank Professor Toru Ohmoto for his instruction to this topic, and Shohei Konno, who was a graduate student at Hokkaido university, for valuable discussions.
The author is also grateful to the editor, associate editor, and reviewers for the helpful comments, some of which suggest potential applications. 
This work was supported by JSPS KAKENHI Grant Number JP22KJ0052, JP25K17299.

\end{document}